\documentclass[reqno, 12pt]{amsart}
\usepackage{amsmath,amsthm,amsfonts,amssymb}
\usepackage[cp1251]{inputenc}
   \usepackage[english]{babel}
\date{}
\newtheorem{theorem}{Theorem}[section]

\newtheorem{corollary}[theorem]{Corollary}

\newtheorem{example}[theorem]{Example}
\numberwithin{equation}{section}

\begin{document}

\centerline{\sc Doubling variables and uniqueness of probability solutions}

\vskip .1in

\centerline{\sc to degenerate stationary Kolmogorov equations}

\vskip .1in

\centerline{\sc V.I. Bogachev, S.V. Shaposhnikov, D.V. Shatilovich}

\vspace*{0.2cm}

{\bf Abstract.} We obtain sufficient conditions for the uniqueness of  a probability solution to the stationary   Kolmogorov equation
with a degenerate  diffusion  matrix. We employ the method of doubling variables known in stochastic analysis
directly to the   Kolmogorov equation.

\vspace*{0.2cm}


\section{Introduction}

The goal of this paper is to study the uniqueness of probability solutions
to the  stationary (elliptic) Kolmogorov equation (also called the Fokker--Planck--Kolmogorov equation)
with a diffusion matrix $A$ and a drift coefficient $b$ in the case where $A$ can be strongly degenerate and is not continuous.

Set
$$
Lf(x)={\rm trace}\bigl(A(x)D^2f(x)\bigr)+\langle b(x), \nabla f(x)\rangle,
$$
where $A(x)=(a^{ij}(x))_{1\le i,j\le d}$ is the  diffusion  matrix, $b(x)=(b^i(x))_{1\le i\le d}$ is the   drift coefficient,
the functions $a^{ij}$ and $b^i$ are Borel  measurable on $\mathbb{R}^d$ and the matrix $A(x)$ is symmetric and nonnegative definite.
Let $d_1\ge 1$ and let
$$\Sigma(x)=(\sigma^{ij}(x))_{1\le i\le d,\, 1\le j\le d_1}$$
be a matrix such that
its  elements $\sigma^{ij}$ are Borel functions on $\mathbb{R}^d$ and the equality
$$
A(x)=\Sigma(x)\Sigma(x)^t
$$
holds. For example, one can take $\Sigma(x)=\sqrt{A(x)}$ and $d_1=d$.
We use the symbol $\langle \cdot,\cdot \rangle$ for the standard inner product on~$\mathbb{R}^d$,
and $|\,\cdot\,|$ denotes the corresponding norm.

A Borel measure $\mu$ on $\mathbb{R}^d$ is called a probability measure if $\mu\ge 0$ and $\mu(\mathbb{R}^d)=1$.
A probability solution to the stationary   Kolmogorov equation
\begin{equation}\label{e1}
L^{*}\mu=0
\end{equation}
 is a  Borel probability measure $\mu$ on $\mathbb{R}^d$ such that the functions $a^{ij}$, $b^i$ are integrable with respect to $\mu$
on every compact set $K\subset\mathbb{R}^d$ and
for every function $f$ from $C_0^{\infty}(\mathbb{R}^d)$ the equality
\begin{equation}\label{int1}
\int_{\mathbb{R}^d}Lf(x)\mu(dx)=0
\end{equation}
holds. Note that this setting allows coefficients that are singular with respect to Lebesgue measure.

An important example of  a probability solution to the  Kolmogorov equation  is an invariant  measure of the diffusion process
with generator $L$. Suppose we are given a homogeneous diffusion process with the family of transition probabilities
$P(t, x, dy)$. The formula
$$
T_tf(x)=\int_{\mathbb{R}^d}f(y)P(t, x, dy)
$$
defines the transition semigroup on the space of bounded Borel functions.
A Borel  probability measure $\mu$ is an invariant measure for this diffusion process
if for every bounded Borel function $f$ on $\mathbb{R}^d$ the equality
$$
\int_{\mathbb{R}^d}T_tf(x)\mu(dx)=\int_{\mathbb{R}^d}f(x)\mu(dx)
$$
holds.

Under broad assumptions about the coefficients $A$ and $b$, the invariant measure $\mu$ satisfies
 equality (\ref{int1}) for all functions $f\in C_0^{\infty}(\mathbb{R}^d)$.
If every probability solution to the  Kolmogorov equation possesses a  continuous
 positive density with respect to  Lebesgue measure, then the corresponding diffusion process
has at  most one invariant probability measure. Indeed, if two probability measures
$\mu_1$ and $\mu_2$ are invariant, then  the measure $|\mu_1-\mu_2|$ is also invariant (see \cite[Lemma 5.1.4]{bookFPK}),
then it satisfies the   Kolmogorov equation
and has a continuous positive density, which is impossible, since the continuous density of the
measure $\mu_1-\mu_2$ must have zeroes, because its integral over $\mathbb{R}^d$ vanishes. In the  case of smooth coefficients
the existence of a continuous positive density of a probability solution is ensured by the nondegeneracy of the
diffusion matrix or by H\"ormander's condition.
Note that in the general case not every probability solution to the  Kolmogorov equation
is an invariant measure for the semigroup, because there are examples (see \cite[Chapter 4]{bookFPK}) in which the Kolmogorov
 equation  with the unit matrix $A$ and a smooth vector field  $b$ has infinitely many linearly independent
probability solutions. If the matrix $A$ is nondegenerate, the functions $a^{ij}$ have Sobolev derivatives,
and the functions $b^i$ are locally Lebesgue integrable to a power higher than the dimension,
then the existence of probability solution to the  Kolmogorov equation
enables us to construct a semigroup and a diffusion process with generator $L$
(see \cite[Chapter 5]{bookFPK}, \cite{Lee-Trutnau}, and~\cite{Lee-Stannat-Trutnau}).

Suppose that the coefficients $a^{ij}$ and $b^i$ are merely Borel  measurable.
Then, according  to the superposition principle  (see \cite[Theorem 4.6]{ZV23}, \cite{BRSsup23}, and \cite{Trev}),
for every probability solution $\mu$ on $\mathbb{R}^d$ satisfying the condition
$$
\int_{\mathbb{R}^d}\frac{\|A(x)\|+|\langle b(x), x\rangle|}{1+|x|^2}\mu(dx)<\infty,
$$
there exists a  Borel probability measure $P$ on $C([0, +\infty), \mathbb{R}^d)$
such that for every $t\ge 0$ and every Borel set $B$ the equality
$P(\omega\colon \omega(t)\in B)=\mu(B)$ holds and $P$ is a solution to the martingale problem with the operator $L$, i.e.,
for every function $f\in C_0^{\infty}(\mathbb{R}^d)$ the process
$$
\xi_t(\omega)=f(\omega(t))-f(\omega(0))-\int_0^tLf(\omega(s), s)\,ds
$$
is a martingale with respect to $P$ and the filtration  $\{\mathcal{F}_t\}=\sigma(\omega(s), s\le t)$.
Moreover, according to  \cite[Chapter 4, Proposition~2.1]{Ikeda-Watanabe}, one can find a probability space
$(\widetilde{\Omega}, \widetilde{\mathcal{F}}, \widetilde{P})$ with a filtration $\{\widetilde{\mathcal{F}}_t\}$, a process $X$
adapted to $\{\widetilde{\mathcal{F}}_t\}$  and an $\widetilde{\mathcal{F}}_t$-Brownian motion $W=(W_t)_{t\ge 0}$
for which $\mu=\widetilde{P}\circ X_t^{-1}$ for all $t\ge 0$ and
$$
dX_t=\sqrt{2}\Sigma(X_t)dW_t+b(X_t)\,dt.
$$
Note that in this generality for some initial conditions the stochastic differential
 equation can have several solutions, and for some initial conditions it can fail to have
solutions. Thus, a probability solution to the  Kolmogorov equation has a probabilistic interpretation in the
 form of  one-dimensional distributions of a solution to the stochastic equation even in the situation, where there is no
semigroup and there is no classical diffusion process defined by means of a family of transition probabilities.

Sufficient conditions for the existence of a probability solution are well known.
In the classical theorem of Hasminskii \cite{H1}, \cite{H2} it is assumed that the
 coefficients are Lipschitz continuous and there exists a Lyapunov function. The condition of Lipschitzness can be substantially
weakened (see \cite[Chapter 2]{bookFPK},  \cite{ZV23}, \cite{Lee2023}),
for example, along with the existence of a Lyapunov function it suffices to have the following:
the matrix $A$ is nondegenerate, its entries belong to the class VMO
and the coefficient $b$ is  integrable with respect to   Lebesgue measure to a sufficiently high power.
The Kolmogorov  equation with a partially degenerate  matrix $A$ is studied in~\cite{ShatShap23}.
In this paper we do not assume that the diffusion matrix is nondegenerate.
In this case it suffices that the coefficients be continuous.  According to \cite[Corollary~2.4.4]{bookFPK}), the  following assertion is true:
if the functions $a^{ij}$ and $b^i$ are continuous on $\mathbb{R}^d$
and there exists a Lyapunov function $V\in C^2(\mathbb{R}^d)$ satisfying the conditions
$\lim\limits_{|x|\to+\infty}V(x)=+\infty$ and $LV(x)\le -C$ whenever $|x|>R$ for some $C>0$ and $R>0$, then the
equation $L^{*}\mu=0$ has a probability solution. There is no uniqueness of a probability solution under such conditions. If $A$ and $b$ vanish
on a nonempty set, then every probability measure with support in this set is a
probability solution.

Let us consider the following example of non-uniqueness of a probability solution to the stationary equation
with a diffusion matrix degenerate at a single point in the presence of a Lyapunov function.

\begin{example}\label{exa1}
\rm
Let $d=3$. Then the equation
$$
\Delta(|x|^2\mu)-{\rm div}(-|x|^2\mu x)=0
$$
admits two different probability solutions $\mu=\delta_0$ and $\mu=C |x|^{-2}e^{-|x|^2/2}\,dx$,
although there is a simple Lyapunov function $V(x)=|x|^2$, for which
$$LV(x)=8|x|^2-2|x|^4\le 16-|x|^4.$$

If we do not assume that the coefficients are smooth, but merely continuous, then replacing $|x|^2$ with $|x|^{\alpha}$ with 
$\alpha<2$ we obtain an example for $d=2$, and for $\alpha<1$ also an example for $d=1$. In addition, the 
drift $b(x)$ up to the positive factor $(2+\alpha)$ equals 
$-\nabla |x|^{2+\alpha}$, i.e., is the gradient of the strictly concave function $-|x|^{2+\alpha}$,
in particular, $\langle b(x)-b(y), x-y\rangle<0$ whenever $x\neq y$.  Note also that in this example any 
probability solution has finite moments of all orders. This follows from standard estimates with the Lyapunov function $|x|^p$, where $p\ge 2$ (see \cite[Theorem 2.3.2]{bookFPK}).
\end{example}

In the case where the coefficients are sufficiently regular and  the matrix $A$ is nondegenerate,
many different sufficient conditions for the  uniqueness of a probability solution are known
(see \cite[Chapter 4]{bookFPK}, \cite{ZV23}, and \cite{LiSt}).
In the case of a degenerate matrix $A$ the uniqueness of a probability solution is less studied.
Moreover, usually the question of uniqueness of
invariant measures was investigated rather than that of a  solution to the   Kolmogorov equation
(see, for example, \cite{Arnold}, \cite{Cannarsa-DaPrato-Frankowska}, \cite{Kliemann}, and \cite{Veretennikov}).
In this  paper we study the uniqueness of a  probability solution to the stationary   Kolmogorov
 equation with a degenerate  matrix~$A$ not assuming the
smoothness of coefficients and not using H\"ormander's conditions. Our approach is based on the method of doubling
 variables described in \cite{Chen}, \cite[Section 3.3]{LebRi}, and \cite[Lemma 8.1.3]{Strook}.
As an illustration we give an adaptation of the reasoning from the
 proof of \cite[Theorem 2.3]{Chen} to the situation we consider.
Note that it is a priori assumed in  \cite{Chen} that the martingale problem with the operator $L$ is well-posed
and that there exists a solution to the martingale problem for the differential operator obtained after doubling the variables.
We do not use these assumptions, but use the procedure of passing to a common probability space known  in stochastic differential equations
 (see \cite[Chapter 4]{Ikeda-Watanabe}).

Let $\mu$ and $\nu$ be probability solutions to  equation~(\ref{e1}). Suppose that the super\-po\-si\-tion  principle holds.
Then we have a probability space $(\Omega^1, F^1, P^1)$ with a filtration $\{\mathcal{F}^1_t\}$, a probability space
$(\Omega^1, F^1, P^1)$ with a filtration $\{\mathcal{F}^1_t\}$,
 a process $X^1=(X^1_t)_{t\ge 0}$ adapted with respect to $\{\mathcal{F}^1_t\}$, an $\mathcal{F}^1_t$-Brownian motion $W^1=(W^1_t)_{t\ge 0}$,
a process $X^2=(X^2_t)_{t\ge0}$ adapted with respect to  $\{\mathcal{F}^2_t\}$, and an $\mathcal{F}^2_t$-Brownian motion $W^2=(W^2_t)_{t\ge 0}$,
for which the equalities $\mu=P^1\circ(X_t^1)^{-1}$ and $\nu=P^2\circ(X_t^2)^{-1}$ hold for all~$t$ and
$$
dX_t^1=\sqrt{2}\Sigma(X_t^1)dW_t^1+b(X_t^1)\,dt, \quad dX_t^2=\sqrt{2}\Sigma(X_t^2)dW_t^2+b(X_t^2)\,dt.
$$
Once can construct (see \cite[Chapter 4]{Ikeda-Watanabe}) a new  probability space $(\Omega, \mathcal{F}, P)$
with a filtration  $\{\mathcal{F}_t\}$ and  processes $X$, $Y$, $W$, adapted with 
respect to~$\{\mathcal{F}_t\}$, such that  the  process $W$ is an $\mathcal{F}_t$-Brownian motion,
the distribution of $(X, W)$ coincides with the distribution of $(X^1, W^1)$, and the distribution of $(Y, W)$ coincides with the distribution of $(X^2, W^2)$,
in particular, one has
\begin{equation}\label{s1}
\left\{\begin{array}{l}
dX_t=\sqrt{2}\Sigma(X_t)dW_t+b(X_t)\,dt, \\
dY_t=\sqrt{2}\Sigma(Y_t)dW_t+b(Y_t)\,dt.
\end{array}\right.
\end{equation}
Let us give a brief description of the construction of such probability space, following
\cite[Chapter 4]{Ikeda-Watanabe}. This construction will be needed below.
Let $\mathcal{P}^1$ and $\mathcal{P}^2$ be the
distributions of $(X^1, W^1)$ and $(X^2, W^2)$ in the space
$C([0, +\infty), \mathbb{R}^d)\times C([0, +\infty), \mathbb{R}^d)$. The projections of the measures $\mathcal{P}^1$ and $\mathcal{P}^2$
onto the second factor equal the  Wiener measure $P^{W}$. Denote by $Q^1_{w}$ and $Q^2_w$ the conditional measures for $\mathcal{P}^1$ and $\mathcal{P}^2$
with respect to~$P^W$.
This means that the integrals against $\mathcal{P}^1$ are evaluated by the formula
$$
\int f(x,y) \, \mathcal{P}^1(dx\, dy)=\int \int f(u,w)\, Q^1_{w}(du)\, P^W(du)
$$
and similarly for $\mathcal{P}^2$.
Then the measure
$$P(dxdydw)=Q^1_{w}(dx)Q^2_{w}(dy)P^{W}(dw)$$
and the set
$$\Omega=C([0, +\infty), \mathbb{R}^d)\times C([0, +\infty), \mathbb{R}^d)\times C([0, +\infty), \mathbb{R}^d)$$
with the natural filtration  define the desired probability space.

System (\ref{s1}) is associated with the differential  operator
\begin{equation}\label{doubL}
\mathbb{L}\psi(x, y)={\rm trace}\bigl(\mathbb{A}(x, y)D^2\psi(x, y)\bigr)+
\langle\mathbb{B}(x, y), \nabla\psi(x, y)\rangle,
\end{equation}
where
$$
\mathbb{A}(x, y)=\left(
                     \begin{array}{cc}
                       A(x) & \Sigma(x)\Sigma(y)^t \\
                       \Sigma(y)\Sigma(x)^t & A(y) \\
                     \end{array}
                   \right),
$$
$$
\mathbb{B}(x, y)=\left(
                   \begin{array}{c}
                     b(x) \\
                     b(y)
                   \end{array}
                 \right).
$$
Set
$$
q(x, y)=\langle x-y, b(x)-b(y)\rangle+{\rm trace}\bigl(\Sigma(x)-\Sigma(y)\bigr)\bigl(\Sigma(x)-\Sigma(y)\bigr)^t.
$$
Observe that
$$
q(x, y)=\mathbb{L}\frac{|x-y|^2}{2}.
$$
Suppose that  for some number $\lambda>0$ and all $(x, y)$ the  condition
$$
q(x, y)\le -\lambda|x-y|^2
$$
is fulfilled and the coefficients $a^{ij}$ and $b^i$ are such that with the aid of It\^o's formula one can  derive the equality
$$
e^{2\lambda t}\mathbb{E}|X_t-Y_t|^2=\mathbb{E}|X_0-Y_0|^2+
\int_0^t\Bigl(2\lambda e^{2\lambda s}\mathbb{E}|X_s-Y_s|^2+2e^{2\lambda s}\mathbb{E}q(X_s, Y_s)\Bigr)\,ds.
$$
Observe that
$$
2\lambda e^{2\lambda s}\mathbb{E}|X_s-Y_s|^2+2e^{2\lambda s}\mathbb{E}q(X_s, Y_s)\le 0.
$$
Then we obtain the inequality
$$
\mathbb{E}|X_t-Y_t|^2\le e^{-2\lambda t}\mathbb{E}|X_0-Y_0|^2.
$$
Recall that the $2$-Kantorovich metric $W_2(\mu, \nu)$ is defined (see, e.g., \cite{B18})
as the infimum of the function
$$
\pi\mapsto \biggl(\int_{\mathbb{R}^d\times\mathbb{R}^d}|x-y|^2\pi(dxdy)\biggr)^{1/2}
$$
on the set of all probability measures $\pi$ on $\mathbb{R}^d\times\mathbb{R}^d$
with projections $\mu$ and $\nu$ on the factors. It is clear that $W_2(\mu, \nu)^2\le \mathbb{E}|X_t-Y_t|^2$.
Letting $t\to+\infty$, we obtain the equality $\mu=\nu$.

The key role in this reasoning is played by the passage from the processes $(X^1, W^1)$ and $(X^2, W^2)$
on different probability spaces to the processes $X$ and $Y$ on a common probability space with a single Wiener process~$W$.
In this  paper we apply this approach directly to solutions of the stationary   Kolmogorov equation:
for given probability solutions $\mu$ and $\nu$ of equation (\ref{e1}) on $\mathbb{R}^d$ we construct
a probability solution $\pi$ to the equation $\mathbb{L}^{*}\pi=0$ on $\mathbb{R}^d\times\mathbb{R}^d$
with projections $\mu$ and $\nu$ on the factors.
Due to the superposition  principle, for the measure $\pi$ we can construct processes $X$, $Y$ and $W$
for which (\ref{s1}) is fulfilled and the one-dimensional distributions of $X_t$ and $Y_t$ are $\mu$ and $\nu$ for each~$t$.

However, the superposition principle does not guarantee that the distribution of
$(X^1, W^1)$ equals the distribution of $(X, W)$ and the distribution of $(X^2, W^2)$ equals the distribution of $(Y, W)$.
In this sense constructing the measure $\pi$ is a simpler task.
Observe that for the  processes $X$ and $Y$ constructed by the measure $\pi$ not only the
one-dimensional distributions $X_t$ and $Y_t$ do not depend on $t$, but also the  distribution of the vector $(X_t, Y_t)$
does not depend on $t$. This circumstance plays the key role in obtaining  sufficient conditions for the uniqueness of
solutions to the  Kolmogorov  equation. The method of constructing a common  probability
space with the aid of conditional measures described above  does not enable us in the general case to find
processes $X$ and $Y$ for which the distribution of the vector $(X_t, Y_t)$ does not  depend on~$t$.
Let $W=(W_t)_{t\ge0}$ be a one-dimensional Wiener process, $X_0$ a random variable with the
Gaussian  distribution $\gamma$ on $\mathbb{R}$ with zero mean  and variance $1/\sqrt{2}$, where $X_0$ and $W$ are independent.
Set $X_t(\omega)=V_t(X_0, W_{\bullet}(\omega))$, where
$$
V_t(y,w)=ye^{-t}w_t-e^{-t}\int_0^te^{s}w_s\,ds.
$$
Then $dX_t=dW_t-X_t\,dt$, i.e.,  $X_t$ is an  Ornstein--Uhlenbeck  process. In this case
the conditional measures $Q_w$ have the form $Q_w=\gamma\circ V(\,\cdot\,, w)^{-1}$.
Let $X^1=X^2=X$ and $W^1=W^2=W$. Then for the measure $P(dxdydw)=Q_w(dx)Q_w(dy)P^W(dw)$ and any Borel sets $C_1$, $C_2$
we have
\begin{multline*}
P\Bigl((x, y, w)\colon (x(0), y(0))\in C_1\times C_2\Bigr)
\\
=\int_{C([0, +\infty), \mathbb{R})}Q_w(x\colon x(0)\in C_1)Q_w(y\colon y(0)\in C_2)P^W(dw)=
\gamma(C_1)\gamma(C_2).
\end{multline*}
Thus, the corresponding processes $X$, $Y$, $W$ satisfy the system of stochastic differential  equations~(\ref{s1})
with $\Sigma=1/\sqrt{2}$ and $b(x)=-x$, and the vector $(X_0, Y_0)$ has distribution $\gamma\otimes\gamma$.
Since $\mathbb{L}^{*}\bigl(\gamma\otimes\gamma\bigr)\neq 0$, the distribution of the vector $(X_t, Y_t)$ depends on $t$.
Thus, not every passage to a common probability  space and not every solution to the martingale problem with the operator $\mathbb{L}$
ensure independence of the distributions $(X_t, Y_t)$ of~$t$.

In this paper, under very weak restrictions on the coefficients, we prove  that the equation $\mathbb{L}^{*}\pi=0$
has a probability solution with given projections, and then apply this result to obtain
sufficient conditions for the uniqueness of a probability solution for a degenerate   Kolmogorov equation.
Note that the method of doubling variables described  above has been successfully applied  for obtaining sharp
conditions for convergence to the invariant measure in the case of partially degenerate Mckean--Vlasov equations
(see, e.g., \cite{Eberle-Guillin-Zimmer}).
Let us note that in the theory of viscosity solutions an important role is played
by the method of doubling variables, which enables one to substantially generalize
the classical maximum principle (see, e.g., \cite{CIL}).
By using this method and estimates of solutions of the equation adjoint to the Kolmogorov equation,
one can obtain results analogous to those obtained with the aid of the stochastic method of doubling variables described above
(see, e.g. \cite{PP13} and~\cite{P24}).
In this paper we apply the stochastic method of doubling variables.

The paper consists of three sections. Section 2 contains the formulations of our main results and some examples.
The proofs are given in Section~3.

\section{Main results}

We use the function introduced above:
$$
q(x, y)=\langle x-y, b(x)-b(y)\rangle+{\rm trace}\bigl(\Sigma(x)-\Sigma(y)\bigr)\bigl(\Sigma(x)-\Sigma(y)\bigr)^t.
$$
For $x\neq y$ set
$$
r(x, y)=\Big|\Bigl(\Sigma(x)-\Sigma(y)\Bigr)^t\frac{x-y}{|x-y|}\Big|^2.
$$

Here is our main  result.

\begin{theorem}\label{th1}
Suppose that either $q(x, y)>0$ for all $x\neq y$ or $q(x, y)<0$ for all $x\neq y$.
Then there exists at most one probability solution~$\mu$ to the equation $L^{*}\mu=0$ for which
$$
\int_{\mathbb{R}^d}\Bigl(\|\Sigma(x)\|^2+|b(x)||x|\Bigr)\mu(dx)<\infty.
$$
\end{theorem}

This theorem can be illustrated by several examples.

\begin{example}\label{ex1}\rm
Suppose that the coefficients $\sigma^{ij}$ are $b^i$
continuous, the functions $\sigma^{ij}$ are bounded and for all $x \in \mathbb R^d$ and some numbers
$m, C\ge 0$ we have
$$
|b(x)|\le C+C|x|^m.
$$
 Suppose also that for some number $\lambda>0$ and all $x, y\in\mathbb{R}^d$ we have
$$
\langle x-y, b(x)-b(y)\rangle+{\rm trace}\bigl(\Sigma(x)-\Sigma(y)\bigr)\bigl(\Sigma(x)-\Sigma(y)\bigr)^t\le -\lambda|x-y|^2,
$$
Then there exists a unique probability solution to the equation $L^{*}\mu=0$.

Substituting  $y=0$ into the estimate, we obtain  the inequality
$$
\langle x, b(x)\rangle\le C_1-C_2|x|^2
$$
with some positive numbers $C_1$ and $C_2$, which yields  that for every $k\ge 2$
 one has $L|x|^k\le N_k-M_k|x|^k$ for all $x$ and some positive numbers $N_k$, $M_k$.
Now \cite[Corollary~2.4.4]{bookFPK} implies the existence of a probability solution $\mu$.
Applying  \cite[Theorem~2.3.2]{bookFPK}, we obtain that $|x|^k\in L^1(\mu)$ for all natural numbers $k\ge 2$.
The uniqueness of a probability solution follows from Theorem~\ref{th1}.
Observe that in this example the uniqueness can be derived without Theorem~\ref{th1} from \cite[Theorem 2.3]{Chen}
or justified by the reasoning presented in the introduction.
\end{example}

Our next example is quite surprising and is not covered by the known results even in the nondegenerate case.

\begin{example}\label{ex2}\rm
Let $b=0$ and let the functions $\sigma^{ij}$ be bounded. If for all $x\neq y$ we have
$$
\Sigma(x)\neq \Sigma(y),
$$
then for $x\neq y$ we have the inequality
$$
q(x, y)={\rm trace}\bigl(\Sigma(x)-\Sigma(y)\bigr)\bigl(\Sigma(x)-\Sigma(y)\bigr)^t>0.
$$
According to Theorem \ref{th1}, there exists at  most one probability solution to the equation $L^{*}\mu=0$.
For example, if at some point $x_0$ we have $\Sigma(x_0)=0$, then in the
considered situation there are no  probability solutions except for Dirac's measure $\delta_{x_0}$.
We emphasize  that no regularity or nondegeneracy of the matrix $A(x)=\Sigma(x)\Sigma(x)^t$ is assumed,
except for the  boundedness of the matrix $\Sigma$.
Note also that in the one-dimensional case this assertion can be easily  derived without Theorem \ref{th1}.
Let $a(x)=\sigma(x)^2$. The equality $(a\mu)''=0$ implies the equality $a\mu=c_1+c_2x$. Since $a\in L^1(\mu)$, we have $a\mu=0$.
Since $\sigma(x)\neq\sigma(y)$ for $x\neq y$, there exists at most one point $x$ at which $a(x)=0$.
Therefore, either there is no probability solution or a probability solution is unique and is Dirac's measure
at the zero of the function $a$. The multidimensional case is much more complicated, because it is impossible to find solutions explicitly.

Under the same assumptions that the functions  $\sigma^{ij}$ are bounded and $\Sigma(x)\neq \Sigma(y)$ whenever $x\neq y$
the uniqueness assertion extends to the case of a  nonzero drift $b$ satisfying the condition
$$
\langle x-y, b(x)-b(y)\rangle\ge 0.
$$
In this case the uniqueness holds in the class of probability measures with respect to which the  function $|b(x)||x|$ is integrable.
\end{example}

Let us give yet another example of conditions on the coefficients that guarantee the inequality $q(x, y)>0$ for $x\neq y$.

\begin{example}\label{ex3}\rm
Suppose that for some number $\lambda>0$ and all $x\neq y$ we have
$$
\langle x-y, b(x)-b(y)\rangle\ge -\lambda|x-y|^2,
$$
$$
{\rm trace}\bigl(\Sigma(x)-\Sigma(y)\bigr)\bigl(\Sigma(x)-\Sigma(y)\bigr)^t>\lambda|x-y|^2.
$$
Then the  condition $q(x, y)>0$ for $x\neq y$ is fulfilled and by Theorem~\ref{th1} there
exists at most one  solution in the class of probability measures with respect to which the
functions $\|\Sigma(x)\|^2$ and $|b(x)||x|$ is integrable.

Suppose, for example, that the matrix $\Sigma$ is diagonal of dimension  $d\times d$. Then our condition on $\Sigma$ takes the form
$$
{\rm trace}\bigl(\Sigma(x)-\Sigma(y)\bigr)\bigl(\Sigma(x)-\Sigma(y)\bigr)^t=\sum_{k=1}^d\bigl(\sigma^{kk}(x)-\sigma^{kk}(y)\bigr)^2
>\lambda|x-y|^2.
$$
If the mapping $f\colon\mathbb{R}^d\to\mathbb{R}^d$ of the form $f(x)=(\sigma^{11}(x), \ldots, \sigma^{dd}(x))$ satisfies the
estimate $\|f(x)-f(y)\|>\sqrt{\lambda}|x-y|$ for $x\neq y$, then our condition on the matrix $\Sigma$ is fulfilled. Such estimate holds  for any
diffeomorphism $f$ with $\|Df^{-1}\|<\lambda^{-1/2}$.
\end{example}

Our next  assertion enables us to relax the condition $q(x, y)<0$ required in Theorem \ref{th1}.

\begin{theorem}\label{th2}
Suppose that there is a nonnegative Borel function $\Lambda$ on $\mathbb{R}^d$ such that for all $x, y\in\mathbb{R}^d$
we have
$q(x, y)\le \bigl(\Lambda(x)+\Lambda(y)\bigr)|x-y|^2$ and for $x\neq y$ the strict  inequality
$q(x, y)<2r(x, y)$ holds. Then there exists at most one probability solution~$\mu$ to the equation $L^{*}\mu=0$ for which
$$
\int_{\mathbb{R}^d}\Biggl(\frac{\|\Sigma(x)\|^2}{1+|x|^2}+\frac{|b(x)|}{1+|x|}+\Lambda(x)\Biggr)\mu(dx)<\infty.
$$
\end{theorem}

\begin{corollary}\label{col1}
Suppose that there is a  nonnegative Borel function $\Lambda$ on $\mathbb{R}^d$ such that
 $\|\Sigma(x)-\Sigma(y)\|\le \bigl(\sqrt{\Lambda(x)}+\sqrt{\Lambda(y)}\bigr)|x-y|$ for all $x, y\in\mathbb{R}^d$,
and if $x\neq y$, then $q(x, y)<2r(x, y)$.
Then there exists at most one probability solution~$\mu$ to the equation $L^{*}\mu=0$ for which
$$
\int_{\mathbb{R}^d}\Biggl(\frac{|b(x)|}{1+|x|}+\Lambda(x)\Biggr)\mu(dx)<\infty.
$$
\end{corollary}
\begin{proof}
Observe that
$$
\|\Sigma(x)\|\le \|\Sigma(0)\|+\bigl(\sqrt{\Lambda(x)}+\sqrt{\Lambda(0)}\bigr)|x|,
$$
$$
r(x, y)\le 2\bigl(\Lambda(x)+\Lambda(y)\bigr)|x-y|^2.
$$
Therefore, $(1+|x|)^{-1}\|\Sigma(x)\|\in L^2(\mu)$ and $q(x, y)\le 4\bigl(\Lambda(x)+\Lambda(y)\bigr)|x-y|^2$,
which enables us to apply Theorem \ref{th2}.
\end{proof}

\begin{example}\label{ex4}\rm
Let $A(x)=\sigma(x)^2I$, where $\sigma$ is a Lipschitz function.
Then
$$
q(x, y)=d\bigl(\sigma(x)-\sigma(y)\bigr)^2+\big\langle x-y, b(x)-b(y)\big\rangle, \quad r(x, y)
=\bigl(\sigma(x)-\sigma(y)\bigr)^2,
$$
and the condition $q(x, y)<2r(x, y)$ takes the form
$$
\langle x-y, b(x)-b(y)\rangle<(2-d)\bigl(\sigma(x)-\sigma(y)\bigr)^2.
$$
Note that for $d=2$ this inequality does not involve the function $\sigma$, hence it imposes no
restrictions on the function $\sigma$ except for its Lipschiztness. Thus, in the case $d=2$,
for every Lipschitz function $\sigma$ and every vector field $b$ satisfying the condition
$$
\langle b(x)-b(y), x-y\rangle<0
$$
for all $x\neq y$, there exists at most one solution in the class of probability measures with respect to which
 the function $|b(x)|(1+|x|)^{-1}$ is integrable.

Example \ref{exa1} shows that the restriction $d=2$ in this observation is sharp:
in all dimensions $d\ge3$ uniqueness fails, and if $\sigma$ is not Lipschitz, then it also fails in dimension $d=2$.

\end{example}

As we already mentioned in the introduction, the main idea of the proof of Theorems \ref{th1} and \ref{th2}
consists in doubling  variables and constructing a probability solution to a new  Kolmogorov equation with the
differential operator $\mathbb{L}$ defined by formula (\ref{doubL}).
Hence the key role in our reasoning is played by the following assertion.

\begin{theorem}\label{main}
Let $\mu$ and $\nu$ be two probability measures on $\mathbb{R}^d$ such that $L^{*}\mu=0$, $L^{*}\nu=0$ and
$$
\int_{\mathbb{R}^d}\Biggl(\frac{\|\Sigma(x)\|^2}{1+|x|^2}+\frac{|b(x)|}{1+|x|}\Biggr)(\mu+\nu)(dx)<\infty.
$$
Then there exists a probability measure $\pi$ on $\mathbb{R}^d\times\mathbb{R}^d$ such that
$\mathbb{L}^{*}\pi=0$ and the projection of the measure $\pi$ on the first factor  equals  $\mu$ and the
 projection on the second factor equals $\nu$.
\end{theorem}

\section{ Proofs of the main results}

\begin{proof}[Proof of Theorem~\ref{main}]
Let $\omega\in C_0^{\infty}(\mathbb{R}^d)$, $\omega\ge 0$, $\omega(x)=0$ whenever $|x|>1$
and $\|\omega\|_{L^1(\mathbb{R}^d)}=1$. For $\varepsilon\in(0, 1)$ we set
$\omega_{\varepsilon}(x)=\varepsilon^{-d}\omega(x/\varepsilon)$. Let $\gamma$ be the standard Gaussian density on~$\mathbb{R}^d$.
Set
$$
\mu_{\varepsilon}(x)=\varepsilon\gamma(x)+(1-\varepsilon)\int_{\mathbb{R}^d}\omega_{\varepsilon}(x-y)\mu(dy),
$$
$$
a^{ij}_{\varepsilon, \mu}(x)=\frac{1-\varepsilon}{\mu_{\varepsilon}(x)}\int_{\mathbb{R}^d}a^{ij}(y)\omega_{\varepsilon}(x-y)\mu(dy)
+\frac{\varepsilon\gamma(x)\delta^{ij}}{\mu_{\varepsilon}(x)},
$$
$$
\sigma^{ij}_{\varepsilon, \mu}(x)=\frac{1-\varepsilon}{\mu_{\varepsilon}(x)}\int_{\mathbb{R}^d}\sigma^{ij}(y)\omega_{\varepsilon}(x-y)\mu(dy),
$$
$$
b^i_{\varepsilon, \mu}(x)=\frac{1-\varepsilon}{\mu_{\varepsilon}(x)}\int_{\mathbb{R}^d}b^{i}(y)\omega_{\varepsilon}(x-y)\mu(dy)
-\frac{\varepsilon\gamma(x)x_i}{\mu_{\varepsilon}(x)}.
$$
Set $A_{\varepsilon, \mu}=(a^{ij}_{\varepsilon, \mu})$,
$\Sigma_{\varepsilon, \mu}=(\sigma^{ij}_{\varepsilon, \mu})$ and $b_{\varepsilon, \mu}=(b_{\varepsilon, \mu}^i)$.
Consider the operator
$$
L_{\varepsilon, \mu}\psi(x)={\rm trace}\bigl(A_{\varepsilon, \mu}(x)D^2\psi(x)\bigr)
+\langle b_{\varepsilon, \mu}(x), \nabla\psi(x)\rangle.
$$
The measure $\mu_{\varepsilon}\, dx$, which will be identified with its density~$\mu_\varepsilon$,  is a solution of the
equation $L_{\varepsilon, \mu}^{*}\mu_{\varepsilon}=0$.

The coefficients of the  operator $L_{\varepsilon, \mu}$ are smooth functions and $\det A_{\varepsilon, \mu}>0$.
Whenever $|x-y|<1$, we have the inequalities
$$
\frac{1}{1+|x|}\le \frac{2}{1+|y|} \quad \hbox{\rm and} \quad \frac{1}{1+|x|^2}\le \frac{3}{1+|y|^2},
$$
which yield the estimate
$$
\frac{|a^{ij}_{\varepsilon, \mu}(x)|}{1+|x|^2}\le
\frac{3(1-\varepsilon)}{\mu_{\varepsilon}(x)}\int_{\mathbb{R}^d}\frac{|a^{ij}(y)|}{1+|y|^2}\omega_{\varepsilon}(x-y)\mu(dy)
+\frac{\varepsilon\gamma(x)}{(1+|x|^2)\mu_{\varepsilon}(x)},
$$
where
$$
\frac{|b^i_{\varepsilon, \mu}(x)|}{1+|x|}=\frac{2(1-\varepsilon)}{\mu_{\varepsilon}(x)}\int_{\mathbb{R}^d}
\frac{|b^{i}(y)|}{1+|y|}\omega_{\varepsilon}(x-y)\mu(dy)
+\frac{\varepsilon\gamma(x)}{\mu_{\varepsilon}(x)}.
$$
Therefore,
$$
\int_{\mathbb{R}^d}\Biggl(\frac{|a^{ij}_{\varepsilon, \mu}(x)|}{1+|x|^2}+\frac{|b^i_{\varepsilon, \mu}(x)|}{1+|x|}\Biggr)\mu_{\varepsilon}(x)\,dx<\infty.
$$
According to \cite[Proposition 5.3.9]{bookFPK}, there exists a Lyapunov  function  $V_{\varepsilon, \mu}\in C^2(\mathbb{R}^d)$
(in the cited proposition the function $V$ belongs to the Sobolev  class $W^{d+, 2}_{loc}$, but in the
case of smooth coefficients one can take a twice continuously differentiable function $V$)
such that $L_{\varepsilon, \mu}V_{\varepsilon, \mu}(x)\to-\infty$ as $|x|\to\infty$.

Similarly we define  $\nu_{\varepsilon}$, $A_{\varepsilon, \nu}=(a^{ij}_{\varepsilon, \nu})$,
$\Sigma_{\varepsilon, \nu}=(\sigma^{ij}_{\varepsilon, \nu})$, $b_{\varepsilon, \nu}=(b^i_{\varepsilon, \nu})$,
and the differential operator $L_{\varepsilon, \nu}$, and for this operator and the measure
$\nu_{\varepsilon}(x)\,dx$ there is also a Lyapunov function $V_{\varepsilon, \nu}$.

Set
$$
\mathbb{L}_{\varepsilon}\psi(x, y)={\rm trace}\bigl(\mathbb{A}_{\varepsilon}(x, y)D^2\psi(x, y)\bigr)+
\langle\mathbb{B}_{\varepsilon}(x, y), \nabla\psi(x, y)\rangle,
$$
where
$$
\mathbb{A}_{\varepsilon}(x, y)=\left(
                     \begin{array}{cc}
                       A_{\varepsilon, \mu}(x) & \Sigma_{\varepsilon, \mu}(x)\Sigma_{\varepsilon, \nu}(y)^t \\
                       \Sigma_{\varepsilon, \nu}(y)\Sigma_{\varepsilon, \mu}(x)^t & A_{\varepsilon, \nu}(y) \\
                     \end{array}
                   \right),
$$
$$
\mathbb{B}(x, y)=\left(
                   \begin{array}{c}
                     b_{\varepsilon, \mu}(x) \\
                     b_{\varepsilon, \nu}(y) \\
                   \end{array}
                 \right).
$$
Let us verify that the matrix $\mathbb{A}_{\varepsilon}$ is nonnegative definite. We have
$$
\Big\langle \mathbb{A}_{\varepsilon}\Big(
                                      \begin{array}{c}
                                        \xi \\
                                        \eta \\
                                      \end{array}
                                    \Big), \Big(
                                      \begin{array}{c}
                                        \xi \\
                                        \eta \\
                                      \end{array}
                                    \Big)
\Big\rangle=\big\langle A_{\varepsilon, \mu}\xi, \xi\big\rangle+
\big\langle \Sigma_{\varepsilon, \nu}^t\eta, \Sigma_{\varepsilon, \mu}^t\xi\big\rangle+
\big\langle \Sigma_{\varepsilon, \mu}^t\xi, \Sigma_{\varepsilon, \nu}^t\eta\big\rangle+
\big\langle A_{\varepsilon, \nu}\eta, \eta\big\rangle.
$$
By the   Cauchy inequality
$$
\big\langle \Sigma_{\varepsilon, \nu}(y)^t\eta, \Sigma_{\varepsilon, \mu}(x)^t\xi\big\rangle+
\big\langle \Sigma_{\varepsilon, \mu}(x)^t\xi, \Sigma_{\varepsilon, \nu}^t(y)\eta\big\rangle
\le \big|\Sigma_{\varepsilon, \mu}(x)^t\xi\big|^2+\big|\Sigma_{\varepsilon, \nu}(y)^t\eta\big|^2.
$$
Applying the   Cauchy inequality once again, we obtain
\begin{multline*}
\big|\Sigma_{\varepsilon, \mu}(x)^t\xi\big|^2=
\frac{(1-\varepsilon)^2}{\mu_{\varepsilon}(x)^2}
\sum_j\biggl(\sum_i\int_{\mathbb{R}^d}\sigma^{ij}(z)\xi_i\omega_{\varepsilon}(x-z)\mu(dz)\biggr)^2
\\
\le
\frac{1-\varepsilon}{\mu_{\varepsilon}(x)}\int_{\mathbb{R}^d}\bigl|\Sigma(z)^t\xi\bigr|^2\omega_{\varepsilon}(x-z)\mu(dz).
\end{multline*}
Similarly we  estimate $\big|\Sigma_{\varepsilon, \nu}(y)^t\eta\big|^2$.
Since $A=\Sigma\Sigma^t$, we have
$$
\bigl|\Sigma(z)^t\xi\bigr|^2=\langle A(z)\xi, \xi\rangle
$$
and the following estimate holds:
\begin{multline*}
\big\langle A_{\varepsilon, \mu}\xi, \xi\big\rangle+\big\langle A_{\varepsilon, \nu}\eta, \eta\big\rangle
\\
\ge \frac{1-\varepsilon}{\mu_{\varepsilon}(x)}\int_{\mathbb{R}^d}\Bigl|\Sigma(z)^t\xi\Bigr|^2\omega_{\varepsilon}(x-z)\mu(dz)+
\frac{1-\varepsilon}{\nu_{\varepsilon}(y)}\int_{\mathbb{R}^d}\bigl|\Sigma(z)^t\eta\bigr|^2\omega_{\varepsilon}(y-z)\nu(dz).
\end{multline*}
Therefore, for all $\xi$ and $\eta$ we have
$$
\Big\langle \mathbb{A}_{\varepsilon}\Big(
                                      \begin{array}{c}
                                        \xi \\
                                        \eta \\
                                      \end{array}
                                    \Big), \Big(
                                      \begin{array}{c}
                                        \xi \\
                                        \eta \\
                                      \end{array}
                                    \Big)\Big\rangle\ge 0.
$$
Thus, the matrix $\mathbb{A}_{\varepsilon}$ is symmetric and nonnegative definite.
We observe that
$$
\mathbb{L}_{\varepsilon}\bigl(V_{\varepsilon, \mu}(x)+V_{\varepsilon, \nu}(y)\bigr)=
L_{\varepsilon, \mu}V_{\varepsilon, \mu}(x)+L_{\varepsilon, \nu}V_{\varepsilon, \nu}(y)\to-\infty
$$
as $|x|^2+|y|^2\to+\infty$. Hence the function
$$
\mathcal{V}_{\varepsilon}(x, y)=V_{\varepsilon, \mu}(x)+V_{\varepsilon, \nu}(y)
$$
is a Lyapunov function  for the operator $\mathbb{L}_{\varepsilon}$.
According to \cite[Corollary 2.4.4]{bookFPK}, there exists a probability solution $\pi_{\varepsilon}$
to the equation $\mathbb{L}_{\varepsilon}^{*}\pi_{\varepsilon}=0$.

Set
$$
\bigl(\mathbb{L}_{\varepsilon}\mathcal{V}_{\varepsilon}\bigr)^{+}
=\max\bigl\{\mathbb{L}_{\varepsilon}\mathcal{V}_{\varepsilon}, 0\bigr\}, \quad
\bigl(\mathbb{L}_{\varepsilon}\mathcal{V}_{\varepsilon}\bigr)^{-}
=\max\bigl\{-\mathbb{L}_{\varepsilon}\mathcal{V}_{\varepsilon}, 0\bigr\}.
$$
Since $\mathbb{L}_{\varepsilon}\mathcal{V}_{\varepsilon}(x, y)\to-\infty$ as $|x|^2+|y|^2\to+\infty$, there
exists a number $C_{\varepsilon}>0$ such that
for all $(x, y)$ the following  inequality is fulfilled:
$$
\bigl(\mathbb{L}_{\varepsilon}\mathcal{V}_{\varepsilon}\bigr)^{+}(x, y)\le C_{\varepsilon}.
$$
We have
$$
\mathbb{L}_{\varepsilon}\ln\bigl(1+\mathcal{V}_{\varepsilon}\bigr)=
\frac{1}{1+\mathcal{V}_{\varepsilon}}\mathbb{L}_{\varepsilon}\mathcal{V}_{\varepsilon}
-\frac{1}{\bigl(1+\mathcal{V}_{\varepsilon}\bigr)^2}\Bigl|\sqrt{\mathbb{A}_{\varepsilon}}\nabla\mathcal{V}_{\varepsilon}\Bigr|^2.
$$
Therefore,
$$
\mathbb{L}_{\varepsilon}\ln\bigl(1+\mathcal{V}_{\varepsilon}\bigr)\le C_{\varepsilon}-
\Bigl(\frac{1}{1+\mathcal{V}_{\varepsilon}}\bigl(\mathbb{L}_{\varepsilon}\mathcal{V}_{\varepsilon}\bigr)^{-}
+\frac{1}{\bigl(1+\mathcal{V}_{\varepsilon}\bigr)^2}\Bigl|\sqrt{\mathbb{A}_{\varepsilon}}\nabla\mathcal{V}_{\varepsilon}\Bigr|^2\Bigr).
$$
Applying  \cite[Theorem 2.3.2]{bookFPK}, we obtain
$$
\int_{\mathbb{R}^d_x\times\mathbb{R}^d_y}\Bigl(\frac{1}{1+\mathcal{V}_{\varepsilon}}\bigl|\mathbb{L}_{\varepsilon}\mathcal{V}_{\varepsilon}\bigr|
+\frac{1}{\bigl(1+\mathcal{V}_{\varepsilon}\bigr)^2}\Bigl|\sqrt{\mathbb{A}_{\varepsilon}}\nabla\mathcal{V}_{\varepsilon}\Bigr|^2\Bigr)
\pi_{\varepsilon}(dxdy)<\infty.
$$
Let $f\in C_0^{\infty}(\mathbb{R})$, $0\le f\le 1$, $f(t)=1$ if $|t|<1$ and $f(t)=0$ if $|t|>2$.
Set
$$
\psi_j(x, y)=f\Bigl(\frac{1}{j}\ln\bigl(1+\mathcal{V}_{\varepsilon}(x, y)\bigr)\Bigr).
$$
Let $\varphi\in C_{0}^{\infty}(\mathbb{R}^d)$.
We have
$$
\mathbb{L}_{\varepsilon}\bigl(\varphi(x)\psi_j(x, y)\bigr)=
\psi_j(x, y)L_{\varepsilon, \mu}\varphi(x)
+\varphi(x)\mathbb{L}_{\varepsilon}\psi_j(x, y)+2\langle\mathbb{A}_{\varepsilon}\nabla\psi_j(x, y), (\nabla\varphi(x), 0)\rangle.
$$
Observe that the coefficients of $L_{\varepsilon, \mu}$ depend only on $x$, therefore, the
function $L_{\varepsilon, \mu}\varphi(x)$ is bounded on $\mathbb{R}^d\times\mathbb{R}^d$.
Since
$$
\mathbb{L}_{\varepsilon}\psi_j=\frac{1}{j}f'\,\mathbb{L}_{\varepsilon}\ln\bigl(1+\mathcal{V}_{\varepsilon}\bigr)+
\frac{f''\,\bigl|\sqrt{\mathbb{A}_{\varepsilon}}\nabla\mathcal{V}_{\varepsilon}\bigr|^2}{j^2\bigl(1+\mathcal{V}_{\varepsilon}\bigr)^2}
$$
and the functions on the right-hand side  are integrable against the measure $\pi_{\varepsilon}$,  as shown above,
we have
$$
\lim_{j\to\infty}\int_{\mathbb{R}^d\times\mathbb{R}^d}\varphi\mathbb{L}_{\varepsilon}\psi_j\,\pi_{\varepsilon}(dxdy)=0.
$$
The following  inequality holds:
$$
\big|\langle\mathbb{A}_{\varepsilon}\nabla\psi_j, (\nabla\varphi, 0)\rangle\big|\le
\langle\mathbb{A}_{\varepsilon}\nabla\psi_j, \nabla\psi_j\rangle^{1/2}
\langle\mathbb{A}_{\varepsilon}(\nabla\varphi, 0), (\nabla\varphi, 0)\rangle^{1/2},
$$
where
$$
\langle\mathbb{A}_{\varepsilon}\nabla\psi_j, \nabla\psi_j\rangle^{1/2}=
\frac{|f'|\, |\sqrt{\mathbb{A}_{\varepsilon}}\mathcal{V}_{\varepsilon}|}{1+\mathcal{V}_{\varepsilon}}, \quad
\langle\mathbb{A}_{\varepsilon}(\nabla\varphi, 0), (\nabla\varphi, 0)\rangle^{1/2}
=\langle A_{\varepsilon, \mu}\nabla\varphi, \nabla\varphi\rangle^{1/2}.
$$
Since the function $\varphi$ has compact support, the function $\langle A_{\varepsilon, \mu}\nabla\varphi, \nabla\varphi\rangle^{1/2}$
is bounded. Therefore,
$$
\lim_{j\to\infty}\int_{\mathbb{R}^d\times\mathbb{R}^d}2\langle\mathbb{A}_{\varepsilon}\nabla\psi_j(x, y), \nabla\varphi(x)\rangle\,\pi_{\varepsilon}(dxdy)=0.
$$
Substituting the function $\varphi(x)\psi_j(x, y)$ into the integral  equality defining the solution $\pi_{\varepsilon}$
 and letting $j\to\infty$, we obtain
$$
\int_{\mathbb{R}^d\times\mathbb{R}^d}L_{\varepsilon, \mu}\varphi(x)\pi_{\varepsilon}(dxdy)=0.
$$
Thus, the projection of $\pi_{\varepsilon}$ onto the first factor is a probability solution
to the  Kolmogorov equation with the operator $L_{\varepsilon, \mu}$. However, according to  \cite[Theorem 4.1.6]{bookFPK},
a probability solution is unique and hence this projection equals $\mu_{\varepsilon}(x)\,dx$.
Similarly we show that the projection of $\pi_{\varepsilon}$ onto the second factor equals $\nu_{\varepsilon}(x)\,dx$.

Observe that as $\varepsilon\to 0$ the measures $\mu_{\varepsilon}(x)\,dx$  converge weakly to the measure $\mu$ and the measures
$\nu_{\varepsilon}(x)\,dx$  converge weakly to the measure $\nu$. Therefore, picking a suitable sequence $\varepsilon_j\to 0$,
we can assume that the measures $\pi_{\varepsilon_j}$  converge weakly to a probability measure $\pi$.
The projection of the  measure $\pi$ on the first factor equals $\mu$ and its projection on the second factor equals $\nu$.
Further for notational simplicity we write $\varepsilon$ in place of $\varepsilon_j$.

 Let us verify that $\mathbb{L}^{*}\pi=0$. To this end we have to justify passage to the limit in the integral identity defining the
solution $\pi_{\varepsilon}$.
Let $\widetilde{a}_{\mu}^{ij}$, $\widetilde{a}_{\nu}^{ij}$,
$\widetilde{\sigma}_{\mu}^{ij}$, $\widetilde{\sigma}_{\nu}^{ij}$,
$\widetilde{b}_{\mu}^i$
and  $\widetilde{b}_{\nu}^i$ be some functions of class $C_0^{\infty}(\mathbb{R}^d)$.
Set
$$
\widetilde{A}_{\mu}=(\widetilde{a}_{\mu}^{ij}), \ \widetilde{A}_{\nu}=(\widetilde{a}_{\nu}^{ij}), \
\widetilde{\Sigma}_{\mu}=(\widetilde{\sigma}_{\mu}^{ij}), \ \widetilde{\Sigma}_{\nu}=(\widetilde{\sigma}_{\nu}^{ij}), \
\widetilde{b}_{\mu}=(\widetilde{b}_{\mu}^{i}), \ \widetilde{b}_{\nu}=(\widetilde{b}_{\nu}^{i}).
$$
The matrix $\widetilde{\mathbb{A}}$ is obtained from the matrix $\mathbb{A}$ by replacing  $A(x)$ with $\widetilde{A}_{\mu}(x)$,
$A(y)$ with $\widetilde{A}_{\nu}(y)$, $\Sigma(x)$ with $\widetilde{\Sigma}_{\mu}(x)$
and $\Sigma(y)$ with $\widetilde{\Sigma}_{\nu}(y)$, respectively.
The vector $\widetilde{\mathbb{B}}$ is obtained from the vector $\mathbb{B}$
by replacing $b(x)$ with $\widetilde{b}_{\mu}(x)$ and $b(y)$ with $\widetilde{b}_{\nu}(y)$.
The corresponding differential operator is denoted by~$\widetilde{\mathbb{L}}$.
The expressions
$$
\widetilde{a}^{ij}(x)_{\varepsilon, \mu}, \quad
\widetilde{a}^{ij}(y)_{\varepsilon, \nu}, \quad
\widetilde{\sigma}^{ij}(x)_{\varepsilon, \mu}, \quad
\widetilde{\sigma}^{ij}(y)_{\varepsilon, \nu}, \quad
\widetilde{b}^{i}(x)_{\varepsilon, \mu}, \quad
\widetilde{a}^{i}(y)_{\varepsilon, \nu}
$$
are defined as above, but  the functions $a^{ij}$ are replaced  with $\widetilde{a}^{ij}_{\mu}$ or $\widetilde{a}^{ij}_{\nu}$,
the functions $\sigma^{ij}$ are replaced with  $\widetilde{\sigma}^{ij}_{\mu}$ or $\widetilde{\sigma}^{ij}_{\nu}$,
and the functions $b^i$ are replaced with $\widetilde{b}^i_{\mu}$ or $\widetilde{b}^i_{\nu}$, respectively.
By  $\widetilde{\mathbb{L}}_{\varepsilon}$ we denote the corresponding differential operator.
Since the  functions $\widetilde{a}^{ij}_{\mu}$, $\widetilde{a}^{ij}_{\nu}$, $\widetilde{\sigma}^{ij}_{\mu}$,
$\widetilde{\sigma}^{ij}_{\nu}$,
$\widetilde{b}^i_{\mu}$ and $\widetilde{b}^i_{\sigma}$ belong to $C_0^{\infty}(\mathbb{R}^d)$,
 for every function $\psi\in C_0^{\infty}(\mathbb{R}^d\times\mathbb{R}^d)$
the functions $\widetilde{\mathbb{L}}_{\varepsilon}\psi$ converge uniformly on $\mathbb{R}^d\times\mathbb{R}^d$
 to the function $\widetilde{\mathbb{L}}\psi$ as $\varepsilon\to 0$.
Observe that due to  weak convergence of $\pi_{\varepsilon}$ to $\pi$ we have
$$
\lim_{\varepsilon\to 0}\int_{\mathbb{R}^d\times\mathbb{R}^d}\widetilde{\mathbb{L}}\psi(x, y)\pi_{\varepsilon}(dxdy)=
\int_{\mathbb{R}^d\times\mathbb{R}^d}\widetilde{\mathbb{L}}\psi(x, y)\pi(dxdy).
$$
Therefore,
$$
\lim_{\varepsilon\to 0}\int_{\mathbb{R}^d\times\mathbb{R}^d}\widetilde{\mathbb{L}}_{\varepsilon}\psi(x, y)\pi_{\varepsilon}(dxdy)=
\int_{\mathbb{R}^d\times\mathbb{R}^d}\widetilde{\mathbb{L}}\psi(x, y)\pi(dxdy).
$$
Let us estimate the expressions
$$
\int_{\mathbb{R}^d\times\mathbb{R}^d}\bigl|\mathbb{L}_{\varepsilon}\psi(x, y)-\widetilde{\mathbb{L}}_{\varepsilon}\psi(x, y)
\bigr|\pi_{\varepsilon}(dxdy),
$$
$$
\int_{\mathbb{R}^d\times\mathbb{R}^d}\bigl|\mathbb{L}\psi(x, y)-\widetilde{\mathbb{L}}\psi(x, y)
\bigr|\pi(dxdy).
$$
Let $B$ be an  open ball in $\mathbb{R}^d$ such that the support of $\psi$ belongs to the set $B\times B$.
Denote by $B'$ the open ball with the same center as $B$ and of radius greater by~$1$.
Observe that
\begin{multline*}
\int_{B\times B}\bigl|b_{\varepsilon, \mu}(x)-\widetilde{b}_{\varepsilon, \mu}(x)\bigr|\pi_{\varepsilon}(dxdy)
\\
\le
\int_{B}\bigl|b_{\varepsilon, \mu}(x)-\widetilde{b}_{\varepsilon, \mu}(x)\bigr|\mu_{\varepsilon}(dx)
\le \int_{B'}\bigl|b(x)-\widetilde{b}_{\mu}(x)\bigr|\mu(dx).
\end{multline*}
Similarly, we have
$$
\int_{B\times B}\bigl|b_{\varepsilon, \nu}(y)-\widetilde{b}_{\varepsilon, \nu}(y)\bigr|\pi_{\varepsilon}(dxdy)
\le \int_{B'}\bigl|b(y)-\widetilde{b}_{\nu}(y)\bigr|\nu(dy),
$$
$$
\int_{B\times B}\bigl|a^{ij}_{\varepsilon, \mu}(x)-\widetilde{a}^{ij}_{\varepsilon, \mu}(x)\bigr|\pi_{\varepsilon}(dxdy)
\le \int_{B'}\bigl|a^{ij}(x)-\widetilde{a}^{ij}_{\mu}(x)\bigr|\mu(dx),
$$
$$
\int_{B\times B}\bigl|a^{ij}_{\varepsilon, \nu}(y)-\widetilde{a}^{ij}_{\varepsilon, \nu}(y)\bigr|\pi_{\varepsilon}(dxdy)
\le \int_{B'}\bigl|a^{ij}(y)-\widetilde{a}_{\nu}^{ij}(y)\bigr|\nu(dy)
$$
By the Cauchy  inequality
\begin{multline*}
\int_{B\times B}\|\Sigma_{\varepsilon, \mu}(x)\Sigma_{\varepsilon, \nu}(y)^t
-\widetilde{\Sigma}_{\varepsilon, \mu}(x)\widetilde{\Sigma}_{\varepsilon, \nu}(y)^t\|\pi_{\varepsilon}(dxdy)
\\
\le \biggl(\int_{B}\|\Sigma_{\varepsilon, \nu}(y)\|^2\nu_{\varepsilon}(y)\,dy\biggr)^{1/2}
\biggl(\int_{B}\|\Sigma_{\varepsilon, \mu}(x)-\widetilde{\Sigma}_{\varepsilon, \mu}(x)\|^2\mu_{\varepsilon}(x)\,dx\biggr)^{1/2}
\\
+\biggl(\int_{B}\|\widetilde{\Sigma}_{\varepsilon, \mu}(x)\|^2\mu_{\varepsilon}(x)\,dy\biggr)^{1/2}
\biggl(\int_{B}\|\Sigma_{\varepsilon, \nu}(y)-\widetilde{\Sigma}_{\varepsilon, \nu}(y)\|^2\nu_{\varepsilon}(y)\,dy\biggr)^{1/2}.
\end{multline*}
Therefore,
\begin{multline*}
\int_{B\times B}\|\Sigma_{\varepsilon, \mu}(x)\Sigma_{\varepsilon, \nu}(y)^t
-\widetilde{\Sigma}_{\varepsilon, \mu}(x)\widetilde{\Sigma}_{\varepsilon, \nu}(y)^t\|\pi_{\varepsilon}(dxdy)
\\
\le C(d,d_1)\|\Sigma\|_{L^2(B', \nu)}\|\Sigma-\widetilde{\Sigma}_{\mu}\|_{L^2(B', \mu)}
\\
+C(d, d_1)\|\Sigma\|_{L^2(B', \mu)}\Bigl(\|\Sigma-\widetilde{\Sigma}_{\nu}\|_{L^2(B', \nu)}
+\|\Sigma-\widetilde{\Sigma}_{\nu}\|_{L^2(B', \nu)}^2\Bigr),
\end{multline*}
where the number $C(d, d_1)$ depends only on $d$ and $d_1$.
Combining the obtained estimates, we arrive at the inequality
\begin{multline*}
\int_{\mathbb{R}^d\times\mathbb{R}^d}
\bigl|\mathbb{L}_{\varepsilon}\psi-\widetilde{\mathbb{L}}_{\varepsilon}\psi\bigr|\pi_{\varepsilon}(dxdy)\le
 C(d, d_1, \psi)\Bigl(\|A-\widetilde{A}_{\mu}\|_{L^1(B', \mu)}+\|A-\widetilde{A}_{\nu}\|_{L^1(B', \nu)}\Bigr)
\\
+C(d, d_1, \psi)\|\Sigma\|_{L^2(B', \nu)}\|\Sigma-\widetilde{\Sigma}_{\mu}\|_{L^2(B', \mu)}
+C(d, d_1, \psi)\|\Sigma\|_{L^2(B', \mu)}\|\Sigma-\widetilde{\Sigma}_{\nu}\|_{L^2(B', \nu)}
\\
+C(d, d_1, \psi)\Bigl(\|b-\widetilde{b}_{\mu}\|_{L^1(B', \mu)}+\|b-\widetilde{b}_{\nu}\|_{L^1(B', \nu)}\Bigr),
\end{multline*}
where the number $C(d, d_1, \psi)$ depends on $d$, $d_1$ and $\psi$, but does not depend on $\varepsilon$.
Clearly, the same estimate is true for the integral of the  function $|\mathbb{L}\psi-\widetilde{\mathbb{L}}\psi|$
against the measure $\pi$. Let $\delta>0$. Due to the   estimates obtained above there exist  functions
$\widetilde{a}_{\mu}^{ij}$, $\widetilde{a}_{\nu}^{ij}$,
$\widetilde{\sigma}_{\mu}^{ij}$, $\widetilde{\sigma}_{\nu}^{ij}$, $\widetilde{b}_{\mu}^i$
and $\widetilde{b}_{\nu}^i$ of class $C_0^{\infty}(\mathbb{R}^d)$ such that
$$
\int_{\mathbb{R}^d\times\mathbb{R}^d}
\bigl|\mathbb{L}_{\varepsilon}\psi-\widetilde{\mathbb{L}}_{\varepsilon}\psi\bigr|\pi_{\varepsilon}(dxdy)<\delta, \quad
\int_{\mathbb{R}^d\times\mathbb{R}^d}
\bigl|\mathbb{L}\psi-\widetilde{\mathbb{L}}\psi\bigr|\pi(dxdy)<\delta.
$$
We deal further with positive  $\varepsilon$ so small that
$$
\biggl|\int_{\mathbb{R}^d\times\mathbb{R}^d}\widetilde{\mathbb{L}}\psi(x, y)\pi_{\varepsilon}(dxdy)-
\int_{\mathbb{R}^d\times\mathbb{R}^d}\widetilde{\mathbb{L}}\psi(x, y)\pi(dxdy)\biggr|<\delta.
$$
Then for such $\varepsilon$ we have the estimate
$$
\biggl|\int_{\mathbb{R}^d\times\mathbb{R}^d}\mathbb{L}_{\varepsilon}\psi(x, y)\pi_{\varepsilon}(dxdy)-
\int_{\mathbb{R}^d\times\mathbb{R}^d}\mathbb{L}\psi(x, y)\pi(dxdy)\biggr|<3\delta.
$$
Since  $\delta$ was arbitrary, we obtain the equality
$$
\int_{\mathbb{R}^d\times\mathbb{R}^d}\mathbb{L}\psi(x, y)\pi(dxdy)=0.
$$
Therefore, the measure $\pi$ is a solution to the equation $\mathbb{L}^{*}\pi=0$.
\end{proof}

\begin{proof}[Proof of Theorem \ref{th1}]
Let $\pi$ be the probability solution to the Kolmogorov equation $\mathbb{L}^{*}\pi=0$ constructed in Theorem~\ref{main}.
Observe that
$$
\mathbb{L}\Bigl(\frac{|x-y|^2}{2}\Bigr)=q(x, y).
$$
Let $f\in C_0^{\infty}(\mathbb{R})$, $0\le f\le 1$, $f(t)=1$ if $|t|<1$ and $f(t)=0$ if $|t|>2$.
Set
$$
\psi_j(x, y)=f\Bigl(\frac{1}{j}\bigl(|x|^2+|y|^2\bigr)\Bigr), \quad E_j=\Bigl\{(x, y)\colon j\le |x|^2+|y|^2\le 2j\Bigr\}.
$$
We have the  equalities
\begin{multline*}
\mathbb{L}\Bigl(\psi_j(x, y)\frac{|x-y|^2}{2}\Bigr)=
\psi_j(x, y)q(x, y)
+\Big\langle\mathbb{A}(x, y)\nabla\psi_j(x, y), \nabla\frac{|x-y|^2}{2}\Big\rangle
\\
+\frac{|x-y|^2}{2}\mathbb{L}\psi_j(x, y),
\end{multline*}
\begin{multline*}
\mathbb{L}\psi_j(x, y)=\frac{2}{j}f'\Bigl(\frac{1}{j}\bigl(|x|^2+|y|^2\bigr)\Bigr)
\Bigl({\rm trace}A(x)+\langle b(x), x\rangle+{\rm trace}A(y)+\langle b(y), y\rangle\Bigr)
\\
+\frac{1}{j^2}f''\Bigl(\frac{1}{j}\bigl(|x|^2+|y|^2\bigr)\Bigr)\Big|\sqrt{\mathbb{A}(x, y)}\nabla(|x|^2+|y|^2)\Big|^2.
\end{multline*}
There is a number $C>0$ independent of  $j$ such that for all $(x, y)$ one has
$$
\Big|\Big\langle\mathbb{A}(x, y)\nabla\psi_j(x, y), \nabla\frac{|x-y|^2}{2}\Big\rangle\Big|\le
CI_{E_j}(x, y)\Bigl(\|\Sigma(x)\|^2+\|\Sigma(y)\|^2\Bigr),
$$
$$
\frac{|x-y|^2}{2}\Big|\mathbb{L}\psi_j(x, y)\Big|\le
CI_{E_j}(x, y)\Bigl(\|\Sigma(x)\|^2+\|\Sigma(y)\|^2+|b(x)||x|+|b(y)||y|\Bigr).
$$
By the  equality
$$
\int_{\mathbb{R}^d\times\mathbb{R}^d}\mathbb{L}\Bigl(\psi_j(x, y)\frac{|x-y|^2}{2}\Bigr)\pi(dxdy)=0
$$
we have
\begin{multline*}
\int_{\mathbb{R}^d\times\mathbb{R}^d}\psi_j(x, y)q(x, y)\pi(dxdy)
\\
=-\int_{\mathbb{R}^d\times\mathbb{R}^d}\Bigl(\Big\langle\mathbb{A}(x, y)\nabla\psi_j(x, y), \nabla\frac{|x-y|^2}{2}\Big\rangle
+\frac{|x-y|^2}{2}\mathbb{L}\psi_j(x, y)\Bigr)\pi(dxdy).
\end{multline*}
Letting $j\to\infty$, we obtain
$$
\int_{\mathbb{R}^d\times\mathbb{R}^d}q(x, y)\pi(dxdy)=0.
$$
By assumption either $q(x, y)>0$ for all $x\neq y$ or $q(x, y)<0$ for all $x\neq y$.
In addition, $q(x, x)=0$ for all $x$. Therefore, $x=y$ for $\pi$-almost all~$(x, y)$.
Let $\varphi\in C_0^{\infty}(\mathbb{R}^d)$.
Then
$$
\int_{\mathbb{R}^d}\varphi(x)\mu(dx)=\int_{\mathbb{R}^d\times\mathbb{R}^d}\varphi(x)\pi(dxdy)
=\int_{\mathbb{R}^d\times\mathbb{R}^d}\varphi(y)\pi(dxdy)=\int_{\mathbb{R}^d}\varphi(y)\nu(dy),
$$
which means that $\mu=\nu$.
\end{proof}

\begin{proof}[Proof of Theorem \ref{th2}.]
Let $S(x, y)=\bigl(\Sigma(x)-\Sigma(y)\bigr)^t$.
Set $H(v)=\ln(1+v/\delta)$, where $v\ge 0$ and $\delta>0$.
Observe that
$$
H'(v)=\frac{1}{\delta+v}, \quad H''(v)=-\bigl(H'(v)\bigr)^2=-\frac{1}{(\delta+v)^2}
$$
Then
$$
\mathbb{L}H\Bigl(\frac{|x-y|^2}{2}\Bigr)=H'\Bigl(\frac{|x-y|^2}{2}\Bigr)q(x, y)
+H''\Bigl(\frac{|x-y|^2}{2}\Bigr)\Big|S(x, y)(x-y)\Big|^2.
$$
Again let $f\in C_0^{\infty}(\mathbb{R})$, $0\le f\le 1$, $f(t)=1$ if $|t|<1$ and $f(t)=0$ if $|t|>2$.
Set
$$
\varphi_j(x, y)=f\Bigl(\frac{1}{j}\ln\bigl(1+|x|^2+|y|^2\bigr)\Bigr), \quad
E_j=\Bigl\{(x, y)\colon j\le\ln(1+|x|^2+|y|^2)\le 2j\Bigr\}.
$$
We have
\begin{multline*}
\mathbb{L}\Big(\varphi_j(x, y)^2H\Big(\frac{|x-y|^2}{2}\Big)\Big)
\\
=\varphi_j(x, y)^2\Big(H'\Bigl(\frac{|x-y|^2}{2}\Big)q(x, y)
+H''\Big(\frac{|x-y|^2}{2}\Bigr)\Big|S(x, y)(x-y)\Big|^2\Big)
\\
+4\varphi_j(x, y)H'\Bigl(\frac{|x-y|^2}{2}\Bigr)\Big\langle\mathbb{A}(x, y)\nabla\varphi_j(x, y), \nabla\frac{|x-y|^2}{2}\Big\rangle+
H\Bigl(\frac{|x-y|^2}{2}\Bigr)\mathbb{L}\varphi_j(x, y)^2.
\end{multline*}
Let $0<\theta<1$. There is a number $C(\theta)>0$ such  that
\begin{multline*}
4\varphi_j(x, y)H'\Bigl(\frac{|x-y|^2}{2}\Bigr)\Big\langle\mathbb{A}(x, y)\nabla\varphi_j(x, y), \nabla\frac{|x-y|^2}{2}\Big\rangle
\\
\le -\theta\varphi_j(x, y)^2H''\Bigl(\frac{|x-y|^2}{2}\Bigr)\Big|S(x, y)(x-y)\Big|^2+
C(\theta)I_{E_j}(x, y)\Bigl(\frac{\|\Sigma(x)\|^2}{1+|x|^2}+\frac{\|\Sigma(y)\|^2}{1+|y|^2}\Bigr).
\end{multline*}
For some number $C>0$ we have the estimate
$$
H\Bigl(\frac{|x-y|^2}{2}\Bigr)\mathbb{L}\varphi_j(x, y)^2\le
CI_{E_j}(x, y)\Bigl(\frac{\|\Sigma(x)\|^2}{1+|x|^2}+\frac{\|\Sigma(y)\|^2}{1+|y|^2}
+\frac{|b(x)|}{1+|x|}+\frac{|b(y)|}{1+|y|}\Bigr).
$$
Set
\begin{multline*}
Q_j(x, y)=C(\theta)I_{E_j}(x, y)\Bigl(\frac{\|\Sigma(x)\|^2}{1+|x|^2}+\frac{\|\Sigma(y)\|^2}{1+|y|^2}\Bigr)
\\
+ CI_{E_j}(x, y)\Bigl(\frac{\|\Sigma(x)\|^2}{1+|x|^2}+\frac{\|\Sigma(y)\|^2}{1+|y|^2}
+\frac{|b(x)|}{1+|x|}+\frac{|b(y)|}{1+|y|}\Bigr).
\end{multline*}
Observe that
$$
\lim_{j\to\infty}\int_{\mathbb{R}^d\times\mathbb{R}^d}Q_j(x, y)\pi(dxdy)=0.
$$
Using the estimates obtained above and taking into account the equalities $H'(v)=(\delta+v)^{-1}$ and $H''(v)=-(\delta+v)^{-2}$,
we obtain the inequality
\begin{multline*}
\mathbb{L}\Big(\varphi_j(x, y)^2H\Bigl(\frac{|x-y|^2}{2}\Bigr)\Big)
\\
\le
\frac{\varphi_j(x, y)^2q(x, y)}{\delta+\frac{|x-y|^2}{2}}-
\frac{(1-\theta)\varphi_j(x, y)^2\Big|S(x, y)(x-y)\Big|^2}{\Big(\delta+\frac{|x-y|^2}{2}\Big)^2}+Q_j(x, y).
\end{multline*}
Integrating against the measure $\pi$, we arrive at the estimate
\begin{multline}\label{new-e}
\int_{\mathbb{R}^d\times\mathbb{R}^d}\frac{(1-\theta)\varphi_j(x, y)^2\Big|S(x, y)(x-y)\Big|^2}{\Bigl(\delta+\frac{|x-y|^2}{2}\Bigr)^2}\,\pi(dxdy)
\\
\le
\int_{\mathbb{R}^d\times\mathbb{R}^d}\frac{\varphi_j(x, y)^2q(x, y)}{\delta+\frac{|x-y|^2}{2}}\pi(dxdy)
+\int_{\mathbb{R}^d\times\mathbb{R}^d}Q_j(x, y)\pi(dxdy).
\end{multline}
Since $q(x, y)\le\bigl(\Lambda(x)+\Lambda(y)\bigr)|x-y|^2$, we have
\begin{multline*}
\int_{\mathbb{R}^d\times\mathbb{R}^d}\frac{(1-\theta)\varphi_j(x, y)^2\big|S(x, y)(x-y)\big|^2}{\Bigl(\delta+\frac{|x-y|^2}{2}\Bigr)^2}\,\pi(dxdy)
\\
\le 2\int_{\mathbb{R}^d\times\mathbb{R}^d}\bigl(\Lambda(x)+\Lambda(y)\bigr)\pi(dxdy)
+\int_{\mathbb{R}^d\times\mathbb{R}^d}Q_j(x, y)\pi(dxdy).
\end{multline*}
Letting $j\to\infty$, we deduce that
$$
\int_{\mathbb{R}^d\times\mathbb{R}^d}
\frac{(1-\theta)\big|S(x, y)(x-y)\big|^2}{\Bigl(1+\frac{|x-y|^2}{2\delta}\Bigr)^2}\,\pi(dxdy)<\infty.
$$
From \eqref{new-e} we obtain
\begin{multline*}
\int_{\mathbb{R}^d\times\mathbb{R}^d}
\frac{\varphi_j(x, y)^2\Big(2|S(x, y)(x-y)|^2-q(x, y)|x-y|^2\Big)}{2\Bigl(\delta+\frac{|x-y|^2}{2}\Bigr)^2}\,\pi(dxdy)
\\
\le
\int_{\mathbb{R}^d\times\mathbb{R}^d}\frac{\delta \varphi_j(x, y)^2q(x, y)}{\Bigl(\delta+\frac{|x-y|^2}{2}\Bigr)^2}\pi(dxdy)
\\
+\theta\int_{\mathbb{R}^d\times\mathbb{R}^d}
\frac{\varphi_j(x, y)^2|S(x, y)(x-y)|^2}{\Bigl(\delta+\frac{|x-y|^2}{2}\Bigr)^2}\,\pi(dxdy)
+\int_{\mathbb{R}^d\times\mathbb{R}^d}Q_j(x, y)\pi(dxdy).
\end{multline*}
 Observe that
$$
2|S(x, y)(x-y)|^2-q(x, y)|x-y|^2=I_{x\neq y}(x, y)\bigl(2r(x, y)-q(x, y)\bigr)|x-y|^2.
$$
Applying again the estimate $q(x, y)\le\bigl(\Lambda(x)+\Lambda(y)\bigr)|x-y|^2$
and letting first $j\to\infty$ and then $\theta\to 0$, we obtain
\begin{multline*}
\int_{\mathbb{R}^d\times\mathbb{R}^d}
\frac{I_{x\neq y}(x, y)\bigl(2r(x, y)-q(x, y)\bigr)|x-y|^2}{2\Bigl(\delta+\frac{|x-y|^2}{2}\Bigr)^2}\,\pi(dxdy)
\\
\le
\int_{\mathbb{R}^d\times\mathbb{R}^d}\frac{\delta|x-y|^2}{\Bigl(\delta+\frac{|x-y|^2}{2}\Bigr)^2}
\bigl(\Lambda(x)+\Lambda(y)\bigr)\pi(dxdy).
\end{multline*}
Since
$$
\lim_{\delta\to 0}\frac{\delta|x-y|^2}{\Bigl(\delta+\frac{|x-y|^2}{2}\Bigr)^2}=0 \quad \hbox{\rm and} \quad
\frac{\delta|x-y|^2}{\Bigl(\delta+\frac{|x-y|^2}{2}\Bigr)^2}\le \frac{1}{2},
$$
we have
$$
\lim_{\delta\to 0}\int_{\mathbb{R}^d\times\mathbb{R}^d}\frac{\delta|x-y|^2}{\Bigl(\delta+\frac{|x-y|^2}{2}\Bigr)^2}
\Bigl(\Lambda(x)+\Lambda(y)\Bigr)\pi(dxdy)=0.
$$
Therefore, letting $\delta\to 0$, we arrive at the inequality
$$
\int_{\mathbb{R}^d\times\mathbb{R}^d}
\frac{I_{x\neq y}(x, y)\bigl(2r(x, y)-q(x, y)\bigr)}{|x-y|^2}\,\pi(dxdy)\le 0,
$$
from which due to the  condition $2r(x, y)-q(x, y)>0$ for $x\neq y$ we conclude  that $x=y$ for $\pi$-almost all $(x, y)$.
As in the proof  of Theorem \ref{th1}, we obtain the equality $\mu=\nu$.
\end{proof}

This research is supported by the Russian Science Foundation Grant N 25-11-00007
at the Lomonosov Moscow State University.

V.I. Bogachev$^{1, 2}$ (vibogach@mail.ru), S.V. Shaposhnikov$^{1, 2}$ (starticle@mail.ru), D.V. Shatilovich$^{1}$ (shatiltop@mail.ru)

1. Moscow Lomonosov State University

2. National Research University ``Higher School of Economics''

\end{document}